\documentclass[a4paper,11pt]{article}
\usepackage{amsmath,amsfonts,amssymb,amsthm,enumerate,titlesec}
\usepackage{fancyhdr}

\newcommand{\ntitle}{On simple ringoids}

\newcommand{\nauthor}{Jens Zumbr\"agel}

\fancypagestyle{plain}{%
  \fancyhead{} 
}

\newcommand{\emphbf}[1]{\textbf{\emph{#1}}}

\newcommand{\N}{\mathbb{N}}
\newcommand{\Z}{\mathbb{Z}}

\newcommand{\on}[1]{\operatorname{#1}}
\newcommand{\id}{\on{id}}
\newcommand{\Aut}{\on{Aut}}
\newcommand{\End}{\on{End}}
\newcommand{\Mult}{\on{Mult}}
\newcommand{\Sym}{\on{Sym}}

\newcommand{\pre}{\preccurlyeq}

\renewcommand{\phi}{\varphi}

\newcommand{\downset}{\:\downarrow\!}

\newtheorem{thm}{Theorem}[section]
\newtheorem{cor}[thm]{Corollary}
\newtheorem{prop}[thm]{Proposition}
\newtheorem{lem}[thm]{Lemma}

\theoremstyle{definition}

\newtheorem{defn}[thm]{Definition}

\theoremstyle{remark}

\newtheorem{exa}[thm]{Example}
\newtheorem{rem}[thm]{Remark}

\newenvironment{enumi}
{\begin{enumerate}[\quad\upshape (i)]%
    \setlength{\itemsep}{1pt}%
    \setlength{\parskip}{0pt}%
    \setlength{\parsep}{0pt}%
}{\end{enumerate}}

\titleformat{\section}
{\large\bfseries}{\thesection.}{.5em}{}

\titleformat{\subsection}
{\large\itshape}{\thesubsection.}{.5em}{}

\def\ccode#1{\par
  \vspace*{8pt}
  {\footnotesize{\leftskip18pt\rightskip\leftskip
      \noindent #1\par}}\par}
\def\abstract#1{\ccode{{\it Abstract}\/.\ #1}}
\def\keywords#1{\ccode{{\it Keywords}\/:\ #1}}

\long\def\symbolfootnote[#1]#2{\begingroup%
\def\thefootnote{\fnsymbol{footnote}}\footnote[#1]{#2}\endgroup} 

\begin{document}

\thispagestyle{empty}

{\centering\Large\bf\ntitle\\}\bigskip

{\centering\large\sc\nauthor%
  \symbolfootnote[1]{School of Mathematical Sciences, University College
    Dublin, Belfield, Dublin 4, Ireland.  E-mail: jens.zumbragel@ucd.ie}%
  \\}\bigskip

\abstract{%
  A ringoid is a set with two binary operations that are linked by the
  distributive laws.  We study special classes of ringoids that are
  congruence-simple or ideal-simple.  In particular, we examine
  generalised parasemifields and non-associative semirings.%
  \symbolfootnote[0]{This work has been supported by the Science
    Foundation Ireland under grant no.\ 08/IN.1/I1950.}}

\keywords{ringoid, semiring, parasemifield, congruence-simple,
  ideal-simple}

\ccode{2000 Mathematics Subject Classification: 16Y99}\medskip


\section{Introduction}

A \emph{ringoid} $(S,+,*)$ is defined as a set $S$ with two binary
operations, $+$ and $*$, such that $a*(b+c)=a*b+a*c$ and
$(a+b)*c=a*c+b*c$ for all $a,b,c\in S$.  In other words, a ringoid is
the most general ring-like structure in which the distributive laws
hold.

In this paper a ringoid $(S,+,*)$ is called a \emph{semiring} if the
groupoid $(S,+)$ is a commutative semigroup; if $(S,+)$ is even an
abelian group we speak of a \emph{ring}.  Semirings or rings $(S,+,*)$
in which the groupoid $(S,*)$ is a semigroup are called
\emph{associative}.  If we have a ringoid $(S,+,*)$ in which $(S,*)$
is a quasigroup we speak of a \emph{generalised parasemifield}; if in
addition $(S,+,*)$ is a semiring the ringoid $(S,+,*)$ is called a
\emph{parasemifield}.

In ringoids the notions of ideal and congruence can be defined (see
Section~\ref{sec_pre}) and we can speak of \emph{ideal-simple} and
\emph{congruence-simple} ringoids.  In the last years there has been a
considerable interest in congruence-simple associative semirings.  The
classification in the commutative case was given in~\cite{ElB01}, and
finite non-commutative semirings with zero were classified
in~\cite{Zum08}.  Studies of general congruence-simple associative
semirings can be found e.g.\ in~\cite{Mon04,ElB07}.  On the other
side, parasemifields appear to be an important class of ideal-simple
semirings and a closer investigation of commutative parasemifields has
been started recently~\cite{Kep08a,Kep08b}.

The aim of the present paper is to initiate the study of more general
ringoids that are congruence-simple or ideal-simple.  We will be
mainly interested in ringoids in which at least one of the binary
operation is non-associative.

The outline of this paper is as follows.  In Section~\ref{sec_pre} we
give the necessary definitions and present some preliminary results.
We present a basic classification of ideal-simple ringoids, provided
that the multiplication is associative and commutative.  In
Section~\ref{sec_gpsf} we examine generalised parasemifields, where in
particular we study the additive groupoid and its automorphism group.
Finally, in Section~\ref{sec_sr} we investigate congruence-simple
ringoids, where we restrict ourselves to the case of (not-necessarily
associative) semirings.  We give a basic classification of
congruence-simple semirings and a characterisation of finite so-called
$k$-ideal-simple semirings having idempotent addition.


\section{Prerequisites}\label{sec_pre}

\subsection*{Groupoids}

\begin{defn}
  A \emphbf{groupoid} $(S,\circ)$ is a non-empty set $S$ with a binary
  operation $\circ:S\times S\to S$.
\end{defn}

We introduce some notation.  Let $\Sym(S)$ be the symmetric group on
the set~$S$, and let $T(S)$ be the monoid of all maps $S\to S$.  The
group of automorphisms of $(S,\circ)$ is denoted by $\Aut(S,\circ)$,
and the monoid of endomorphisms of $(S,\circ)$ is denoted by
$\End(S,\circ)$.  Clearly, $\Aut(S,\circ)$ is a subgroup of $\Sym(S)$
and $\End(S,\circ)$ is a submonoid of $T(X)$.

For $a\in S$ consider the left and right multiplication maps $L_a$ and
$R_a$, defined by $L_a(x)=a\circ x$ and $R_a(x)=x\circ a$ for $x\in
S$.  The submonoid of $T(X)$ generated by $\{L_a\mid a\in S\} \cup
\{R_a\mid a\in S\}$ is denoted by $\Mult(S,\circ)$.

For subsets $A$ and $B$ of $S$ let $A\circ B$ denote the subset
$\{a\circ b\mid a\in A,b\in B\}$ of $S$.  Furthermore, for $a,b\in S$
we will write $a\circ B$ and $A\circ b$ for the sets $\{a\}\circ B$
and $A\circ\{b\}$.

\begin{defn}
  Let $(S,\circ)$ be a groupoid.  A \emphbf{congruence} of $(S,\circ)$
  is an equivalence relation $\rho$ on the set $S$ such that
  $\rho(a,b)$ implies $\rho(x\circ a,x\circ b)$ and $\rho(a\circ
  x,b\circ x)$ for all $a,b,x\in S$.

  A non-empty subset $A$ of $S$ is called a \emphbf{subgroupoid} of
  the groupoid $(S,\circ )$ if $A\circ A\subseteq A$ holds.  It is
  called an \emphbf{ideal} of the groupoid $(S,\circ )$ if we have
  $(S\circ A)\cup(A\circ S)\subseteq A$.
\end{defn}

A subgroupoid or ideal $A$ of $(S,\circ)$ will be called \emph{proper}
if $A\neq S$.

\subsection*{Ringoids}

\begin{defn}
  A \emphbf{ringoid} $(S,+,*)$ is a non-empty set $S$ with two binary
  operations, $+$ and $*$, such that $a*(b+c)=a*b+a*c$ and
  $(a+b)*c=a*c+b*c$ hold for all $a,b,c\in S$.
\end{defn}

As usual we call the operation $+$ \emph{addition} and the operation
$*$ \emph{multiplication}.

\begin{rem} 
  Given two groupoids $(S,+)$ and $(S,*)$, then $(S,+,*)$ is a ringoid
  if and only if $\Mult(S,*)$ is a submonoid of $\End(S,+)$.
\end{rem}

\subsection*{Simple ringoids}

\begin{defn}
  Let $(S,+,*)$ be a ringoid.  A \emphbf{congruence} on $(S,+,*)$ is
  an equivalence relation on the set $S$ that is a congruence of both
  the groupoid $(S,+)$ and the groupoid $(S,*)$.

  A non-empty subset $A$ of $S$ is called an \emphbf{ideal} of the
  ringoid $(S,+,*)$ if $A$ is a subgroupoid of the groupoid $(S,+)$
  and an ideal of the groupoid $(S,*)$; that is, it holds $(A+A)\cup
  (S*A)\cup (A*S)\subseteq A$.
\end{defn}

\begin{defn}
  The ringoid $(S,+,*)$ is said to be \emphbf{congruence-simple} if
  its only congruences are $\rho=\id_S:=\{(a,a)\mid a\in S\}$ and
  $\rho=S\times S$.  It is said to be \emphbf{ideal-simple} if there
  are no proper ideals $A$ of $(S,+,*)$ with $|A|\geq 2$.  If there
  are no proper ideals of $(S,+,*)$ at all, the ringoid will be called
  \emphbf{ideal-free}.
\end{defn}

We note that a ringoid $(S,+,*)$ is congruence-simple if and only if
every non-constant homomorphism $(S,+,*)\to (T,+,*)$ into another
ringoid $(T,+,*)$ is injective.

\subsection*{Ideal-simple ringoids}

We present a basic classification of ideal-simple ringoids, provided
that the multiplication is associative and commutative (note that
these structures are formally akin to semirings, with the role of
addition and multiplication interchanged).

We will use the following standard result:

\begin{lem}
  Let $(S,*)$ be a semigroup such that $S*a=S=a*S$ for all $a\in S$.
  Then $S$ is a group.
\end{lem}

\begin{prop}
  Let $(S,+,*)$ be a ringoid with associative and commutative
  multiplication.  If $(S,+,*)$ is ideal-simple then one of the
  following holds:
  \begin{enumi}
    \item $|S*S|=1$,
    \item $(S,*)$ is a group,
    \item there is an absorbing element $o$ in $(S,*)$ such that
      $(S\setminus\{o\},*)$ is a group.
    \end{enumi}
\end{prop}

\begin{proof}
  We may assume $|S|\geq 2$.  Consider the set \[N:=\{a\in S\mid
  |S*a|=1\}\:,\] which is either empty or an ideal of $(S,+,*)$.  If
  $N$ is empty then for all $a\in S$ we have $|S*a|>1$, and, since
  $S*a$ is an ideal of $(S,+,*)$, we have $S*a=S$.  Consequently,
  $(S,*)$ is a group.  If $N=S$ we have $|S*a|=1$ for all $a\in S$ and
  it follows that $|S*S|=1$.

  Suppose now that $N$ is non-empty and $N\neq S$.  Then $|N|=1$, so
  $N=\{o\}$ for some $o\in S$, and we have $S*o=\{c\}$ for some $c\in
  S$.  Furthermore, for all $a\neq o$ we have $S*a=S$.  Let
  $T:=S\setminus\{o\}$ and $a,b\in T$.  Then $S=S*b=S*a*b$, so that
  $a*b\neq o$.  Since $S*a=S$ and $S*o=\{c\}$ we have $S\setminus\{c\}
  \subseteq T*a \subseteq T$.  This proves that $c=o$ and $T*a=T$ for
  all $a\in T$.  Hence $o$ is absorbing in $(S,*)$ and $(T,*)$ is a
  group.
\end{proof}

In the next section we examine ringoids $(S,+,*)$ in which the
multiplication is a group or, more generally, a quasigroup in more
detail.


\section{Generalised parasemifields}\label{sec_gpsf}

A groupoid $(S,*)$ is called a \emph{quasigroup} if for all $a\in S$
the multiplication maps $L_a:x\mapsto a*x$ and $R_a:x\mapsto x*a$ are
bijective maps of $S$, i.e.\ $L_a,R_a\in\Sym(S)$.  
In this case $\Mult(S,*)$ is a submonoid of $\Sym(S)$.

\begin{defn}
  A \emphbf{generalised parasemifield} is a ringoid $(S,+,*)$ such
  that $(S,*)$ is a quasigroup.
\end{defn}

Apparently, every generalised parasemifield is an ideal-free ringoid.

\begin{rem}\label{psf_rem}
  Given two groupoids $(S,+)$ and $(S,*)$, then $(S,+,*)$ is a
  generalised parasemifield if and only if $\Mult(S,*)$ is a submonoid
  of $\Aut(S,+)$.
\end{rem}

\subsection*{The additive automorphism group}

\begin{lem}\label{psf_lem}
  Let $(S,+,*)$ be a generalised parasemifield.  Then the additive
  automorphism group $\Aut(S,+)$ is transitive.
\end{lem}


\begin{proof}
  Let $x,y\in S$ be given.  Since $R_x\in\Sym(S)$ there exists $a\in
  S$ such that $R_x(a)=y$.  Now $L_a\in\Aut(S,+)$ and we have
  $L_a(x)=a*x=R_x(a)=y$.  
\end{proof}

\begin{cor}\label{psf_cor1}
  Let $(S,+,*)$ be a generalised parasemifield, $|S|>1$.  Then $(S,+)$
  has neither a neutral element nor an absorbing element.
\end{cor}

\begin{proof}
  If $(S,+)$ had a neutral resp.\ absorbing element then, since
  $\Aut(S,+)$ is transitive, every element of $S$ would be neutral
  resp.\ absorbing in $(S,+)$.  But each groupoid can have at most one
  neutral or absorbing element.
\end{proof}

\begin{cor}
  Let $(S,+,*)$ be a generalised parasemifield such that $(S,+)$ is a
  commutative semigroup.  Then $|S|=\infty$.
\end{cor}

\begin{proof}
  Suppose $|S|<\infty$.  It is well-known that every finite semigroup
  has an idempotent element, hence $(S,+)$ has an idempotent element.
  Since $\Aut(S,+)$ is transitive, every element of $S$ is idempotent.
  It follows that $(S,+)$ has an absorbing element, namely $\sum_{x\in
    S}x$, in contradiction to Corollary~\ref{psf_cor1}.
\end{proof}

\pagebreak

\begin{exa}
  The ringoid $(\Z,\max,+)$ is an example of an
  infinite generalised parasemifield with associative and commutative
  addition.
\end{exa}

\subsection*{Groupoids with full automorphism group}

The discussion of the automorphism group of groupoids is motivated by
Remark~\ref{psf_rem} and Lemma~\ref{psf_lem}.

If for a groupoid $(S,+)$ every bijection of $S$ is an automorphism,
i.e.\ $\Aut(S,+) = \Sym(S)$, then for any quasigroup $(S,*)$ we get a
generalised parasemifield $(S,+,*)$.  A proof of the following result
can be found in~\cite{Tam63,Tam67}.

\begin{prop}
  A groupoid $(S,\circ)$ satisfies $\Aut(S,\circ) = \Sym(S)$ if and
  only if $(S,\circ)$ is either isomorphic or anti-isomorphic to one
  of the following types:
  \begin{enumi}
  \item A right zero semigroup, i.e.\ $x\circ y=y$ for
    all $x,y\in S$.  
  \item The idempotent quasigroup of order 3.  
  \item The groupoid $(\{1,2\},\circ)$ where $x\circ 1=2$ and $x\circ
    2=1$, for $x\in\{1,2\}$.
  \end{enumi}
  Furthermore, if $\Aut(S,\circ)$ is triply transitive then
  $\Aut(S,\circ) = \Sym(S)$.
\end{prop}

\subsection*{Groupoids with transitive automorphism group}

We give some remarks on groupoids $(S,\circ)$ with transitive
automorphism group that make computer searches for such groupoids more
efficient.

Fix a groupoid $(S,\circ)$.  For $s\in S$ define the sets
\begin{eqnarray*}
&N_l(s) := \{x\in S\mid s\circ x = x\}\:,\quad
 N_r(s) := \{x\in S\mid x\circ s = x\}\:,\\
&A_l(s) := \{x\in S\mid s\circ x = s\}\:,\quad
 A_r(s) := \{x\in S\mid x\circ s = s\}\:.
\end{eqnarray*}
These are the sets of all elements $x$ for which $s$ is left-neutral,
right-neutral, left-absorbing and right-absorbing, respectively.
Further let
\begin{eqnarray*}
&n_l(s):=|N_l(s)|\:,\quad
 n_r(s):=|N_r(s)|\:,\\
&a_l(s):=|A_l(s)|\:,\quad
 a_r(s):=|A_r(s)|\:.
\end{eqnarray*}

\begin{lem}\label{trans_lem}
  Let $(S,\circ)$ be a groupoid with transitive automorphism group.
  Then the cardinalities $n_l(s)$, $n_r(s)$, $a_l(s)$ and $a_r(s)$ are
  independent of $s$.
\end{lem}

\begin{proof}
  Let $s,t\in S$, $s\neq t$.  Then there exists
  $\phi\in\Aut(S,\circ)$ such that $\phi(s)=t$.  Now, for any
  $x\in S$ we have $s\circ x = x$ if and only if $\phi(s)\circ
  \phi(x) = \phi(x)$, hence $x\in N_l(s)$ if and only if
  $\phi(x)\in N_l(\phi(s)) = N_l(t)$.  This shows $\phi(N_l(s)) =
  N_l(t)$ and thus $n_l(s)=n_l(t)$.

  Similarly one observes $\phi(N_r(s))=N_r(t)$, $\phi(A_l(s))=A_l(t)$
  and $\phi(A_r(s))=A_r(t)$, hence $n_r(s)=n_r(t)$, $a_l(s)=a_l(t)$
  and $a_r(s)=a_r(t)$.
\end{proof}

In the situation of Lemma~\ref{trans_lem} we may just write $n_l$,
$n_r$, $a_l$ and $a_r$ for the corresponding cardinalities.

\begin{lem}
  Let $(S,\circ)$ be a finite groupoid with transitive automorphism
  group.  Then $n_l=a_r$ and $a_l=n_r$.
\end{lem}

\begin{proof}
  For any groupoid $(S,\circ)$ it holds
  \[ \sum_{x\in S}n_l(x) = |\{(x,y)\in S^2\mid x\circ y=y\}| =
  \sum_{y\in S}a_r(y)\:. \] If $(S,\circ)$ has a transitive
  automorphism group it follows $|S|n_l = |S|a_r$ by
  Lemma~\ref{trans_lem}.  Hence, if $|S|$ is finite we have $n_l=a_r$.
  The proof of $a_l=n_r$ is very similar.
\end{proof}

\begin{cor}
  Let $(S,\circ)$ be a commutative finite groupoid with transitive
  automorphism group.  Then $n_l=n_r=a_l=a_r$.
\end{cor}

\begin{proof}
  We have $N_l=N_r$ and $A_l=A_r$ and thus $n_l=n_r$ and $a_l=a_r$.
\end{proof}

\begin{exa}
  Let $m$ be an odd positive integer and let $(\Z_m,+,\cdot)$ be the
  ring of integers modulo $m$.  Define an operation $\circ$ on the set
  $\Z_m$ by
  \[ (a\circ b) := \frac{1}{2}\cdot(a+b) \:,\] where $\frac{1}{2}$ is the
  inverse of $2$ in $\Z_m$.  

  Then $(\Z_m,\circ)$ is a commutative quasigroup with transitive
  automorphism group (and $n_l=n_r=a_l=a_r=1$).  In fact, for every
  $c\in\Z_m$ the translation map $\Z_m\to\Z_m$, $x\mapsto x+c$ is an
  automorphism of $(\Z_m,\circ)$.  We get a generalised parasemifield
  $(\Z_m,\circ,+)$ where both operations are commutative, the addition
  is a quasigroup and the multiplication is a group.
\end{exa}



\section{Congruence-simple semirings}\label{sec_sr}

In this section congruence-simple semirings in which the
multiplication is not necessarily associative will be studied.

\begin{defn}
  A \emphbf{semiring} is a ringoid $(S,+,*)$ in which the additive
  groupoid $(S,+)$ is associative and commutative.
\end{defn}

\subsection*{Basic classification}

Let $(S,+,*)$ be a semiring.  For $a\in S$ and
$n\in\N=\{1,2,3,\dots\}$ we define $na$ recursively by $1a:=a$ and
$(n+1)a:=na+a$.

\begin{rem}
  Let $(S,+,*)$ be a semiring.  For any $a,b\in S$ and $m,n\in\N$ we
  have $m(na) = (mn)a$ and $n(a+b)=na+nb$, which follows from the
  associativity and commutativity of the addition.  Furthermore, it
  holds $a*(nb) = n(a*b)$ and $(na)*b = n(a*b)$, which is implied by
  the distributive laws.
\end{rem}

\begin{lem}\label{sem_lem}
  Let $(S,+,*)$ be a semiring.  Let the relation $\pre$ on $S$ be
  defined by \[a\preccurlyeq b\quad :\Leftrightarrow\quad
  \exists\,n\in\N, x\in S: nb = x+a\] and let
  $\rho\,:=\,\preccurlyeq\cap\preccurlyeq^{-1}$.  Then $\rho$ is a
  congruence of the semiring $(S,+,*)$ and $\rho(a,2a)$ for all
  $a\in S$.
\end{lem}

\begin{proof}
  First we show that $\pre$ is a preorder on $S$.  Let $a,b,c\in S$.
  Since $2a=a+a$ we have $a\pre a$.  Now suppose $a\pre b$ and $b\pre
  c$.  Then there exist $m,n\in\N$ and $x,y\in S$ such that $nb=x+a$
  and $mc=y+b$.  If follows $(nm)c = n(mc) = n(y+b) = ny+nb = ny+x+a$,
  so that $a\pre c$.

  Next we show that the preorder $\pre$ is compatible with the
  semiring operations.  Let $a,b\in S$ such that $a\pre b$, and let
  $s\in S$.  There exist $n\in\N$ and $x\in S$ such that $nb=x+a$.  It
  follows $(n+1)(s+b) = ns+nb+s+b = ns+x+a+s+b = (ns+x+b)+s+a$, so
  that $s+a\pre s+b$.  Furthermore, $n(s*b) = s*(nb) = s*(x+a) =
  s*x+s*a$ and $n(b*s) = (nb)*s = (x+a)*s = x*s+a*s$, so that $s*a\pre
  s*b$ and $a*s\pre b*s$.

  From what we have shown it follows that $\rho=\{(a,b)\in S^2\mid
  a\pre b\text{ and }b\pre a\}$ is a congruence of the semiring
  $(S,+,*)$.  Let $a\in S$.  Then $a\pre 2a$, since $1(2a)=a+a$, and
  $2a\pre a$, since $3a=a+2a$, hence $\rho(a,2a)$.
\end{proof}

\begin{cor}
  Let $(S,+,*)$ be a congruence-simple semiring with an additively
  neutral element $0$.  Then the groupoid $(S,+)$ is either idempotent
  or a group.
\end{cor}

\begin{proof}
  Consider the congruence $\rho$ of Lemma~\ref{sem_lem}.  If
  $\rho=\id_S$, then $2a=a$ for all $a\in S$, i.e.\ $(S,+)$ is
  idempotent.  If $\rho\neq\id_S$ then by congruence-simplicity it
  follows that $\rho=S\times S$, hence for all $a,b\in S$ there exists
  $n\in\N$ and $x\in S$ such that $nb=x+a$.  In particular, for all
  $a\in S$ there exists $x\in S$ such that $0=n0=x+a$.  Thus $(S,+)$
  is a group.
\end{proof}

\begin{rem}
  Suppose that $(S,+,*)$ is a congruence-simple semiring without an
  additively neutral element such that $(S,+)$ is not idempotent.
  Along the same lines as the proof of Theorem~3.1 in~\cite{ElB01}
  one can show that either the groupoid $(S,+)$ is cancellative or
  there is an absorbing element $o$ of $(S,+)$ such that $2x=o$ for
  all $x\in S$.
\end{rem}

\subsection*{Ideals and congruences}

We examine the relationship between ideals, congruences and $k$-ideals
(defined below).  Afterwards we give a classification of $k$-ideal-simple
semirings.

\begin{lem}
  Let $A$ be an ideal of a semiring $(S,+,*)$.  Define a relation
  $\rho_A$ on~$S$ by
  \[\rho_A := \{(x,y)\in S^2\mid \exists\, a,b\in A: x+a=y+b\}\:.\]
  Then $\rho_A$ is a congruence of the semiring $(S,+,*)$.
\end{lem}

\begin{proof}
  The relation $\rho_A$ is clearly reflexive and symmetric.  Let
  $x,y,z\in S$ and suppose $\rho_A(x,y)$ and $\rho_A(y,z)$.  Hence,
  there are $a,b,c,d\in A$ such that $x+a=y+b$ and $y+c=z+d$.  It
  follows $x+a+c = y+b+c = y+c+b = z+d+b$.  Since $a+c,b+d\in A$ we
  have $\rho_A(x,z)$ and thus $\rho_A$ is an equivalence relation.

  To show that $\rho_A$ is a congruence, let $x,y\in S$ such that
  $\rho_A(x,y)$, and let $a,b\in A$ such that $x+a=y+b$.  Further, let
  $z\in S$.  Then $z+x+a = z+y+b$ and thus $\rho_A(z+x,z+y)$.  Also,
  $z*x+z*a = z*y+z*b$ and $x*z+a*z = y*z+b*z$ with $z*a,z*b,a*z,b*z\in
  A$.  Thus $\rho_A(z*x,z*y)$ and $\rho_A(x*z,y*z)$.
\end{proof}

\begin{defn}
  Let $(S,+,*)$ be a ringoid.  A non-empty subset $A$ of $S$ is called
  a \emphbf{$k$-ideal} of $(S,+,*)$ if $A$ is an ideal of the ringoid
  $(S,+,*)$ such that $(A+A^c)\cup (A^c+A)\subseteq A^c$.  Here, $A^c$
  denotes the complement $S\setminus A$ of $A$.

  A ringoid $(S,+,*)$ is called \emphbf{$k$-ideal-simple} if it has no
  proper $k$-ideals $A$ with $|A|\geq 2$.
\end{defn}

\begin{lem}
  Let $(S,+,*)$ be a semiring with an additively neutral element $0$
  and let $A$ be a $k$-ideal.  Then the $\rho_A$-class $[0]_{\rho_A}$
  of $0$ equals $A$.
\end{lem}

\begin{proof}
  We have $[0]_{\rho_A} = \{x\in S\mid \exists\, a\in A: x+a\in A\}$.
  This set equals~$A$, since $a\in A$ and $x+a\in A$ implies $x\in A$
  by the $k$-ideal property.
\end{proof}

\begin{cor}
  Let $(S,+,*)$ be a semiring such that $(S,+)$ has a neutral element.
  If the semiring is congruence-simple then it is $k$-ideal simple.
\end{cor}

\subsection*{Idempotent $k$-ideal-simple semirings}

We consider now semirings with idempotent addition.

\begin{rem}
  For any commutative idempotent semigroup $(S,+)$ there is an order
  relation $\leq$ on $S$ defined by $\leq\;:=\{(a,b)\in S^2\mid
  a+b=b\}$.  The partially ordered set $(S,\leq)$ is then a
  join-semilattice, where $a\vee b=a+b$ for all $a,b\in S$.

  If $(S,+,*)$ is a semiring with idempotent addition and $\leq$ is
  defined as before, we have that $a\leq b$ implies $x*a\leq x*b$ and
  $a*x\leq b*x$ for any $a,b,x\in S$.
\end{rem}
  
The following proposition gives an easily checkable criterion
characterising finite additively idempotent $k$-ideal-simple
semirings, which can be used to accelerate computer searches for
congruence-simple semirings.

\begin{prop}
  Let $(S,+,*)$ be a finite semiring with idempotent addition and let
  $\infty:=\sum_{a\in S}a\in S$.  Let $M$ be the set of minimal
  elements in $(S,\leq)$.  Then $S$ is $k$-ideal simple if and only if
  for all $x\in S\setminus(M\cup\{\infty\})$ we have $\infty*x\not\leq
    x$ or $x*\infty\not\leq x$.
\end{prop}

\begin{proof}
  First we show that a non-empty subset $A$ of $S$ satisfies
  \[A+A\subseteq A\quad\text{ and }\quad A+A^c\subseteq A^c\] if and
  only if $A$ is of the form $A_x := \downset x := \{a\in S\mid a\leq
  x\}$ for some $x\in S$.

  Indeed, if $A+A\subseteq A$ and $A+A^c\subseteq A^c$ we let
  $x:=\sum_{a\in A}a\in A$.  We have $A = \downset x$, since each
  $a\in A$ satisfies $a\leq x$, and conversely, if $a\leq x$ then
  $a+x=x\in A$ and thus $a\in A$.

  On the other hand, let $x\in S$ and $A = \downset x$.  Then
  $A+A\subseteq A$, since $a,b\leq x$ implies $a\vee b\leq x$, for all
  $a,b\in S$.  Furthermore, if $a\in A$ and $y\in S$ with $a+y\in A$
  then $y\leq a+y\in A$, hence $y\in A$.  Thus $A+A^c\subseteq A^c$.

  Next for $x\in S$ we note that $S*A_x\subseteq A_x$ is equivalent to
  $\infty*x\leq x$, since $s*a\leq s*x\leq \infty*x$ for all $s\in S$
  and $a\in A_x$.  Similarly, $A_x*S\subseteq A_x$ if and only if
  $x*\infty\leq x$.  Hence, $A_x$ is a $k$-ideal if and only if
  $\infty*x\leq x$ and $x*\infty\leq x$.  

  We remark that $|A_x|=1$ if and only if $x$ is minimal in
  $(S,\leq)$, i.e.\ $x\in M$.  In conclusion, there exists a proper
  $k$-ideal $A$ with $|A|>1$ if and only if there is $x\in
  S\setminus(M\cup\{\infty\})$ that satisfies the conditions
  $\infty*x\leq x$ and $x*\infty\leq x$.
\end{proof}

\begin{exa}
  There are (up to isomorphism) five simple congruence-semirings
  $(\{0,1,2\},+,*)$ of order 3 with idempotent addition and an
  \emph{absorbing zero} 0, i.e.\ the element $0$ is additively neutral
  and multiplicatively absorbing:

  \[{\small\begin{tabular}{c|ccc}
      + & 0 & 1 & 2 \\
      \hline
      0 & 0 & 1 & 2 \\ 
      1 & 1 & 1 & 2 \\ 
      2 & 2 & 2 & 2 \\ 
    \end{tabular}\qquad	
    \begin{tabular}{c|ccc}
      $*$ & 0 & 1 & 2 \\
      \hline
      0 & 0 & 0 & 0 \\ 
      1 & 0 & 0 & 0 \\ 
      2 & 0 & 2 & 2 \\ 
    \end{tabular}\quad	
    \begin{tabular}{c|ccc}
      $*$ & 0 & 1 & 2 \\
      \hline
      0 & 0 & 0 & 0 \\ 
      1 & 0 & 0 & 1 \\ 
      2 & 0 & 2 & 2 \\ 
    \end{tabular}}\]
  \[{\small\quad	
    \begin{tabular}{c|ccc}
      $*$ & 0 & 1 & 2 \\
      \hline
      0 & 0 & 0 & 0 \\ 
      1 & 0 & 0 & 2 \\ 
      2 & 0 & 0 & 2 \\ 
    \end{tabular}\quad	
    \begin{tabular}{c|ccc}
      $*$ & 0 & 1 & 2 \\
      \hline
      0 & 0 & 0 & 0 \\ 
      1 & 0 & 0 & 2 \\ 
      2 & 0 & 1 & 2 \\ 
    \end{tabular}\quad	
    \begin{tabular}{c|ccc}
      $*$ & 0 & 1 & 2 \\
      \hline
      0 & 0 & 0 & 0 \\ 
      1 & 0 & 0 & 2 \\ 
      2 & 0 & 2 & 2 \\ 
    \end{tabular}}\]  
\end{exa}  

The following table shows the number of congruence-simple semirings
(up to isomorphism) of order $n$ with idempotent addition and an
absorbing zero:

\[\begin{tabular}{c|ccccc}
  $n$ & 2 & 3 & 4 & 5 & 6 \\\hline
  $\#$\,general & 2 & 5 & 428 & 138 167 & ? \\
  $\#$\,commutative & 2 & 1 & 21 & 715 & 59 640 \\
  $\#$\,associative & 2 & 0 & 0 & 0 & 1
\end{tabular}\]


The figures for associative semirings are taken from the
classification~\cite{Zum08}, the other numbers are results of a
computer search program.


\bibliographystyle{amsalpha}
\bibliography{literature}

\end{document}